\documentclass[11pt,twoside]{article}
\usepackage{a4}
\usepackage{amssymb,amsmath,amsthm,latexsym}
\usepackage{amsfonts}
\usepackage[numbers,sort&compress]{natbib}
\usepackage{hyperref}
\usepackage{cleveref}
\usepackage{enumitem}
\usepackage{booktabs}
\usepackage{graphicx}
\usepackage{xcolor}
\usepackage{tikz}
\usetikzlibrary{arrows.meta, positioning, fit, backgrounds}
\usepackage{listings}
\usepackage{float}

\makeatletter
\newcommand{\lean}[1]{%
  \begingroup
    \ttfamily
    \edef\lean@str{\detokenize{#1}}%
    \expandafter\lean@loop\lean@str\lean@stop
  \endgroup
}
\def\lean@stop{\lean@stop}
\def\lean@loop#1{%
  \ifx#1\lean@stop
    \let\lean@next\relax
  \else
    \ifnum`#1>64 \ifnum`#1<91 \allowbreak\fi\fi
    #1%
    \ifnum`#1=95 \allowbreak\fi
    \ifnum`#1=46 \allowbreak\fi
    \let\lean@next\lean@loop
  \fi
  \lean@next
}
\makeatother

\lstdefinelanguage{Lean}{
  keywords={theorem, lemma, def, instance, where, by, have, let, show,
            fun, match, with, if, then, else, do, return, import, open,
            namespace, end, section, variable, noncomputable, sorry,
            structure, class, extends, deriving, inductive, abbrev,
            example, private, protected, local, attribute, set_option,
            obtain, exact, apply, intro, rw, simp, omega, ring, norm_num,
            calc, suffices, constructor, cases, induction, rcases},
  sensitive=true,
  morecomment=[l]{--},
  morecomment=[s]{/-}{-/},
  morestring=[b]",
  literate={→}{$\to\;$}2 {←}{$\leftarrow$}1 {↔}{$\leftrightarrow$}1
           {∀}{$\forall$}1 {∃}{$\exists$}1 {λ}{$\lambda$}1
           {≤}{$\leq$}1 {≥}{$\geq$}1 {≠}{$\neq$}1
           {∈}{$\in$}1 {∉}{$\notin$}1 {⊆}{$\subseteq$}1
           {ℕ}{$\mathbb{N}$}1 {ℤ}{$\mathbb{Z}$}1 {ℚ}{$\mathbb{Q}$}1
           {⟨}{$\langle$}1 {⟩}{$\rangle$}1
           {×}{$\times$}1 {⊤}{$\top$}1 {⊥}{$\bot$}1
           {∧}{$\wedge$}1 {∨}{$\vee$}1 {¬}{$\lnot$}1
           {𝓞}{$\mathcal{O}$}1
           {θ'}{$\theta'$}2 {θ}{$\theta$}1
           {Rˣ}{$R^{\times}$}2
           {↦}{$\mapsto$}1
           {∨}{$\lor\;$}2
           {·}{$\cdot$}1
           {₁}{${}_{1}$}1 {₂}{${}_{2}$}1
}
\newtheorem*{ramanujanquestion}{Question 464}

\newcommand{\Z}{\mathbb{Z}}
\newcommand{\Q}{\mathbb{Q}}
\newcommand{\N}{\mathbb{N}}
\newcommand{\cO}{\mathcal{O}}
\newcommand{\Lean}{\textsf{Lean\,4}}
\newcommand{\Mathlib}{\textsf{Mathlib}}

\begin{document}

\setcounter{page}{1}

\markboth {\hspace*{-9mm} \centerline{\footnotesize \sc
   A formal proof of the Ramanujan--Nagell theorem in \Lean{}}
                 }
                { \centerline{\footnotesize \sc
         Barinder S.\ Banwait} \hspace*{-9mm}
               }

\vspace*{-2cm}

\begin{center}
{
       {\Large \textbf { \sc  A formal proof of the Ramanujan--Nagell theorem in \Lean{}
                               }
       }
\\

\medskip

{\sc Barinder S.\ Banwait}\\
{\footnotesize London, UK}\\
{\footnotesize e-mail: {\it barinder.s.banwait@gmail.com}}
}
\end{center}

\thispagestyle{empty}

\hrulefill

\begin{abstract}
{\footnotesize
We present a complete formalization, in the \Lean{} interactive theorem prover
with the \Mathlib{} library, of the Ramanujan--Nagell theorem: the only
integer solutions to the Diophantine equation $x^2 + 7 = 2^n$
are $(n,x) \in \{(3,\pm1),(4,\pm3),(5,\pm5),(7,\pm11),(15,\pm181)\}$. The formalization
includes all dependencies, notably a proof that the quadratic ring $R = \Z\!\left[(1+\sqrt{-7})/2\right]$ is a Euclidean domain---hence a unique factorization domain---together with the determination of its unit group. The development comprises about 1{,}800 lines of \Lean{} code, contains no \texttt{sorry} placeholders, and depends only on \Mathlib{}'s standard axioms. We describe the proof strategy, the architecture of the formalization, and the challenges
encountered in bridging the gap between textbook proofs and their
machine-checked counterparts, with particular attention to the algebraic
number theory infrastructure required.
}
\end{abstract}
\hrulefill



\section{Introduction}\label{sec:intro}

In 1913, Ramanujan posed the following question to the \emph{Journal of the Indian Mathematical Society}:

\begin{ramanujanquestion}[\cite{Ramanujan1913}]
   $2^n- 7$ is a perfect square for the values $3$, $4$, $5$, $7$, $15$ of $n$. Find other values.
\end{ramanujanquestion} 

Implicit here is the claim that these are the only such values of $n$. This was first proved by Nagell in 1948~\cite{Nagell1948}, and for this reason this result is known in the literature as the \textbf{Ramanujan--Nagell theorem}, often appearing in a first undergraduate course in algebraic number theory.

In this paper we present a complete formalization of the
Ramanujan--Nagell theorem in \Lean{} using the \Mathlib{} library.
The formalization, together with an interactive dependency graph and
proof blueprint, is publicly available at
\begin{center}
\url{https://barindersbanwait.com/ramanujan_nagell/}.
\end{center}
Our contributions are:
\begin{enumerate}[label=(\roman*)]
  \item a machine-checked proof of the theorem itself;
  \item \Lean{} infrastructure for the ring of integers
        $R = \Z[\tfrac{1+\sqrt{-7}}{2}]$ of $\Q(\sqrt{-7})$: a proof that $R$
        is a Euclidean domain---hence a unique factorization domain---and the
        determination of its unit group;
  \item an analysis of the gap between the textbook proof and our formal proof, identifying where informal arguments in the standard mathematical proof required substantial elaboration in Lean.
\end{enumerate}
To our knowledge, this is the first formalization of a conjecture of
Ramanujan in a proof assistant, and the first complete resolution of a
specific \emph{exponential} Diophantine equation in any formal
verification system. Diophantine equations themselves are not new to
formalization: \Mathlib{} already contains a development of the solution
theory of Pell's equation $x^2 - d\,y^2 = 1$, contributed by Michael
Gei{\ss}er and Michael Stoll, growing out of Gei{\ss}er's BSc
thesis~\cite{StollPell}.
What sets the Ramanujan--Nagell equation apart is that the unknown $n$
occurs in an exponent, which calls for quite different methods---here,
the arithmetic of the ring of integers of $\Q(\sqrt{-7})$.

The remainder of this paper is organized as follows.
\Cref{sec:math} recalls the standard mathematical proof.
\Cref{sec:architecture} describes the architecture of the formalization.
\Cref{sec:challenges} discusses the principal challenges encountered.
\Cref{sec:lessons} draws lessons for formalizing algebraic number
theory, including the role of AI-assisted tools.

\section{The mathematical proof}\label{sec:math}

Our exposition follows Stewart and Tall~\cite[Section~4.9]{StewartTall2015}, with
one departure: where they use that $\Q(\sqrt{-7})$ has class number one---so that its
ring of integers is a unique factorization domain---we instead show that
$\Q(\sqrt{-7})$ is norm-Euclidean and obtain unique factorization directly from the
Euclidean algorithm (\Cref{sec:ufd}). We do so because the norm-Euclidean route is
considerably easier to formalize; we compare the two approaches in \Cref{sec:earlier}.

\subsection{Setup and notation}

We work throughout in the ring
\[
  R = \Z\!\left[\tfrac{1+\sqrt{-7}}{2}\right] = \Z \oplus \Z\,\theta,
  \qquad \theta = \frac{1+\sqrt{-7}}{2},
\]
the ring of integers of $\Q(\sqrt{-7})$ (since $-7 \equiv 1 \pmod 4$, the half-integer
$\theta$ is an algebraic integer). Write $\theta' = \frac{1-\sqrt{-7}}{2} = 1 - \theta$
for the conjugate of $\theta$. Directly from these definitions,
\[
  \theta^2 = \theta - 2, \qquad \theta\theta' = 2, \qquad
  \theta + \theta' = 1, \qquad \theta - \theta' = \sqrt{-7}.
\]
For $z = x + y\theta \in R$ with $x, y \in \Z$, the norm $N(z) = z\bar z$
(where $\bar z = x + y\theta'$) is
\[
  N(x + y\theta) = x^2 + xy + 2y^2 = \tfrac14\bigl((2x+y)^2 + 7y^2\bigr).
\]
The second form, obtained by completing the square, makes the nonnegativity
of $N$ manifest---it is a sum of two squares with positive weights---and is the
shape in which $N$ recurs in the unit and Euclidean-domain arguments of
\Cref{sec:ufd}. In particular $N$ is a multiplicative, nonnegative integer that
vanishes only at $z = 0$.

Any solution of $x^2 + 7 = 2^n$ has $x$ odd, since $2^n - 7$ is odd. Moreover $x$ and
$-x$ are solutions for the same $n$, so we assume without loss of generality that
$x > 0$, restoring the sign only in the final tally.

\subsection{The even case}\label{sec:even-case}

If $n$ is even, then $x^2 + 7 = 2^n$ gives a factorization of integers
\[
  (2^{n/2} + x)(2^{n/2} - x) = 7.
\]
Since $7$ is prime, we must have $2^{n/2} + x = 7$ and $2^{n/2} - x = 1$,
giving $2^{1 + n/2} = 8$, hence $n = 4$ and $x = 3$.

\subsection{Unique factorization in $R$}\label{sec:ufd}

The odd case rests on unique factorization in $R$, which we obtain from the Euclidean
algorithm: $\Q(\sqrt{-7})$ is one of the five imaginary quadratic fields that are
\emph{norm-Euclidean}, meaning that $N$ is a Euclidean function on $R$.

To verify this, let $a, b \in R$ with $b \neq 0$ and put
$\gamma = a/b \in \Q(\sqrt{-7})$. It suffices to find $q \in R$ with
$N(\gamma - q) < 1$, for then $N(a - bq) = N(b)\,N(\gamma - q) < N(b)$. Write
$\gamma = \xi + \eta\sqrt{-7}$ with $\xi, \eta \in \Q$; since conjugation sends
$\sqrt{-7} \mapsto -\sqrt{-7}$, any $q = \xi_0 + \eta_0\sqrt{-7}$ has
\[
  N(\gamma - q) = (\xi - \xi_0)^2 + 7\,(\eta - \eta_0)^2.
\]
The weight $7$ on the second coordinate is what makes the rounding delicate: taking
$\xi_0, \eta_0$ to be the nearest \emph{integers} to $\xi, \eta$ gives only
$N(\gamma - q) \le \tfrac14 + 7\cdot\tfrac14 = 2$, which need not be below~$1$. The
half-integers in $R$ supply the missing room. Writing each element of
$R = \Z[\tfrac{1+\sqrt{-7}}{2}]$ as $\tfrac12(j + k\sqrt{-7})$ with $j, k \in \Z$ of
equal parity, its $\sqrt{-7}$-coordinate $k/2$ ranges over $\tfrac12\Z$. So choose
$k \in \Z$ nearest to $2\eta$, giving $|\eta - k/2| \le \tfrac14$, and then $j \in \Z$ of
the same parity as $k$ nearest to $2\xi$, giving $|\xi - j/2| \le \tfrac12$ (the integers
of a fixed parity being spaced $2$ apart). With $q = \tfrac12(j + k\sqrt{-7}) \in R$,
\[
  N(\gamma - q) = \Bigl(\xi - \tfrac{j}{2}\Bigr)^2 + 7\Bigl(\eta - \tfrac{k}{2}\Bigr)^2
    \le \tfrac14 + 7\cdot\tfrac{1}{16} = \tfrac{11}{16} < 1.
\]
Hence $N$ is a Euclidean function on $R$, so $R$ is a Euclidean domain, and therefore
a principal ideal domain (PID) and a unique factorization domain (UFD).

The units of $R$ are precisely $\pm 1$: a unit has norm $1$, and
$\tfrac14\bigl((2x+y)^2 + 7y^2\bigr) = 1$ forces $y = 0$ and $x = \pm 1$. Finally,
$N(\theta) = N(\theta') = 2$ is prime in $\Z$, so $\theta$ and $\theta'$ are
irreducible, hence prime; and they are not associates, since their ratio is not a
unit. We shall use the factorization $2 = \theta\theta'$ repeatedly.

\subsection{Factorization in $R$}\label{sec:factorization}

Now let $n$ be odd with $n \geq 5$, and put $m = n - 2 \geq 3$. Since $x$ is odd,
write $x = 2k+1$; then $x^2 + 7 = 4(k^2 + k + 2)$, so
\begin{equation}\label{eq:divided}
  \frac{x^2 + 7}{4} = k^2 + k + 2 = 2^m.
\end{equation}
Because $x$ is odd, the conjugate elements
\[
  \alpha = \frac{x + \sqrt{-7}}{2} = k + \theta, \qquad
  \beta  = \frac{x - \sqrt{-7}}{2} = k + \theta'
\]
both lie in $R$. Using $\theta\theta' = 2$, their product is
$\alpha\beta = k^2 + k + 2 = 2^m = (\theta\theta')^m = \theta^m\theta'^m$, while their
difference is
\[
  \alpha - \beta = \theta - \theta' = \sqrt{-7}.
\]
The factors $\alpha$ and $\beta$ are coprime in $R$: the only primes dividing
$2 = \theta\theta'$ are $\theta$ and $\theta'$, and a common prime factor of
$\alpha, \beta$ would divide $\alpha - \beta = \sqrt{-7}$; but $N(\sqrt{-7}) = 7$ is
not divisible by $N(\theta) = N(\theta') = 2$. Since $R$ is a unique factorization
domain with units $\pm 1$, comparing the two factorizations of $\theta^m\theta'^m$
forces
\begin{equation}\label{eq:unit-comparison}
  \frac{x + \sqrt{-7}}{2} = \pm\, \theta^m \quad\text{or}\quad \pm\, \theta'^m.
\end{equation}

\subsection{The sign argument}\label{sec:sign}

Regardless of which of the four cases in~\eqref{eq:unit-comparison} holds, by taking the Galois conjugate and subtracting, we obtain the following in every case:
\[
  \pm\sqrt{-7} = \theta^m - \theta'^m.
\]
We claim that the positive sign cannot occur. For, if it does, then we have $\theta^m - \theta'^m = \theta - \theta'$, and working modulo $\theta'^2$ and using $\theta\theta' = 2$ and $\theta = 1 - \theta'$, we have
$\theta^2 = (1 - \theta')^2 \equiv 1 \pmod{\theta'^2}$,
so $\theta^m \equiv \theta \pmod{\theta'^2}$ (since $m$ is odd). Hence
\[
  \theta - \theta' \equiv \theta^m - \theta'^m \equiv \theta \pmod{\theta'^2},
\]
which gives $\theta' \equiv 0 \pmod{\theta'^2}$, a contradiction. Hence the sign is
negative, and we have
\begin{equation}\label{eq:negative-sign}
  -\sqrt{-7} = \theta^m - \theta'^m.
\end{equation}

\subsection{The binomial expansion and reduction modulo~$7$}\label{sec:mod7}

Expanding $\theta^m$ and $\theta'^m$ by the binomial theorem, the
even-indexed terms cancel on subtraction and
\[
  \theta^m - \theta'^m
    = 2^{1-m}\sqrt{-7}\,\Bigl(\binom{m}{1} - \binom{m}{3}\cdot 7
      + \binom{m}{5}\cdot 7^2 - \cdots \pm \binom{m}{m}\cdot 7^{(m-1)/2}\Bigr).
\]
Comparing with the negative-sign identity~\eqref{eq:negative-sign}, cancelling
$\sqrt{-7}$ and multiplying by $2^{m-1}$ gives
\[
  -2^{m-1}
  = \binom{m}{1} - \binom{m}{3}\cdot 7 + \binom{m}{5}\cdot 7^2
    - \cdots \pm \binom{m}{m}\cdot 7^{(m-1)/2},
\]
and reducing modulo~$7$ gives
\begin{equation}\label{eq:mod7}
  -2^{m-1} \equiv m \pmod{7}.
\end{equation}
Since $2^6 \equiv 1 \pmod{7}$, the residue of $-2^{m-1}$ modulo~$7$
depends only on $m \bmod 6$.  A direct check shows that the only
solutions to~\eqref{eq:mod7} modulo~$42$ are
\begin{equation}\label{eq:mod42}
  m \equiv 3,\; 5,\; \text{or } 13 \pmod{42}.
\end{equation}

\subsection{Uniqueness per residue class}\label{sec:uniqueness-math}

It remains to show that each class in~\eqref{eq:mod42} contains at most one value
of~$m$. Suppose, for contradiction, that two solutions give exponents $m < m'$ with
$m \equiv m' \pmod{42}$. Put $d = m' - m > 0$; then $42 \mid d$, so in particular
$7 \mid d$, and write $l = v_7(d) \geq 1$ for the exact power of~$7$ dividing~$d$. All congruences below are modulo~$7^{l+1}$.

We begin with two congruences. First, expanding $(1 + \sqrt{-7})^{d}$ by the
binomial theorem and grouping the terms into a rational part and a $\sqrt{-7}$-part,
\[
  (1 + \sqrt{-7})^{d}
    = \Bigl(1 + \sum_{j \geq 1}\binom{d}{2j}(-7)^{j}\Bigr)
    + \Bigl(d + \sum_{j \geq 1}\binom{d}{2j+1}(-7)^{j}\Bigr)\sqrt{-7}.
\]
Every term with $j \geq 1$ is divisible by~$7^{l+1}$: it carries a factor~$7^{j}$
from $(-7)^{j}$, and from the identity $k\binom{d}{k} = d\binom{d-1}{k-1}$ together
with $7^{l} \mid d$ we get $v_7\binom{d}{k} \geq l - v_7(k)$, so the term has $7$-adic
valuation at least $j + l - v_7(k) \geq l+1$ (since $k \leq 2j+1 < 7^{j}$ forces
$v_7(k) \leq j-1$). Only the two leading terms survive, giving
\begin{equation}\label{eq:binom-cong}
  (1 + \sqrt{-7})^{d} \equiv 1 + d\sqrt{-7} \pmod{7^{l+1}}.
\end{equation}
Second, since $42 \mid d$ we have $3 \mid d$, so $2^{d} - 1 = 8^{d/3} - 1$. The
\emph{lifting-the-exponent lemma}---the identity
$v_7(a^{k} - b^{k}) = v_7(a - b) + v_7(k)$, valid whenever $7 \mid a - b$ but
$7 \nmid ab$---applied with $a = 8$ and $b = 1$ gives
$v_7(2^{d} - 1) = v_7(7) + v_7(d/3) = 1 + l$, so
\begin{equation}\label{eq:two-cong}
  2^{d} \equiv 1 \pmod{7^{l+1}}.
\end{equation}

Dividing~\eqref{eq:binom-cong} by~$2^{d}$ -- permissible because~\eqref{eq:two-cong}
makes $2^{d}$ invertible modulo $7^{l+1}$ -- and recalling $\theta = (1+\sqrt{-7})/2$
gives $\theta^{d} = (1+\sqrt{-7})^{d}/2^{d} \equiv 1 + d\sqrt{-7}$. Applying the
conjugation $\sqrt{-7} \mapsto -\sqrt{-7}$ yields the companion congruence, so
\begin{equation}\label{eq:theta-power-cong}
  \theta^{d} \equiv 1 + d\sqrt{-7}, \qquad
  \theta'^{\,d} \equiv 1 - d\sqrt{-7} \pmod{7^{l+1}}.
\end{equation}

Now both $m$ and $m'$ satisfy the negative-sign identity~\eqref{eq:negative-sign}, so
$\theta^m - \theta'^m = \theta^{m'} - \theta'^{m'}$. Writing $m' = m + d$ and
substituting the congruences~\eqref{eq:theta-power-cong},
\[
  \theta^{m'} - \theta'^{m'}
    = \theta^m\theta^{d} - \theta'^m\theta'^{\,d}
    \equiv \theta^m(1 + d\sqrt{-7}) - \theta'^m(1 - d\sqrt{-7})
    = (\theta^m - \theta'^m) + d\sqrt{-7}\,(\theta^m + \theta'^m).
\]
Since the left-hand side equals $\theta^m - \theta'^m$, the trace
$P = \theta^m + \theta'^m$ satisfies
\begin{equation}\label{eq:trace-cong}
  d\sqrt{-7}\,P \equiv 0 \pmod{7^{l+1}}.
\end{equation}

Finally, $P$ is coprime to~$7$. Indeed the integers $P_k = \theta^k + \theta'^k$
satisfy the recurrence $P_{k+2} = P_{k+1} - 2P_k$ (as $\theta$ and $\theta'$ are the
roots of $t^2 - t + 2$), with $P_0 = 2$ and $P_1 = 1$; modulo~$7$ this sequence is
periodic of period~$3$, cycling through $2, 1, 4$, none of which is~$0$. So $P$ is
invertible modulo~$7^{l+1}$, and cancelling it from~\eqref{eq:trace-cong} leaves
$d\sqrt{-7} \equiv 0 \pmod{7^{l+1}}$. But $\sqrt{-7} = 2\theta - 1$, so
$d\sqrt{-7} = -d + 2d\,\theta$ has coordinates $(-d, 2d)$ in $\Z \oplus \Z\theta$;
divisibility by the rational integer~$7^{l+1}$ thus forces $7^{l+1} \mid d$,
contradicting $v_7(d) = l$. Hence no two solutions share a residue class
in~\eqref{eq:mod42}, and each class contains at most one value of~$m$.

\subsection{Conclusion of the proof}

By~\eqref{eq:mod42} and the uniqueness just proved, each of the residue classes
$m \equiv 3, 5, 13 \pmod{42}$ contributes at most one value of~$m$ from a solution,
and these are realised by $m = 3, 5, 13$, that is $n = 5, 7, 15$, with
$x = 5, 11, 181$ respectively (as one checks directly). Together with the remaining
odd case $n = 3$ (where $x^2 = 1$) and the even case $n = 4$ (\Cref{sec:even-case},
where $x^2 = 9$), this accounts for every solution with $x > 0$. Restoring the sign
of~$x$, the complete set of integer solutions of $x^2 + 7 = 2^n$ is
\[
  (n, x) \in \{(3, \pm 1),\ (4, \pm 3),\ (5, \pm 5),\ (7, \pm 11),\ (15, \pm 181)\}.
\]

\section{Architecture of the formalization}\label{sec:architecture}

\subsection{Overview}\label{sec:overview}

The main theorem is stated in \texttt{Basic.lean} as follows:

\begin{lstlisting}
theorem RamanujanNagell :
  ∀ x : ℤ, ∀ n : ℕ, x ^ 2 + 7 = 2 ^ n →
    (x, n) = (1, 3) ∨ (x, n) = (-1, 3)
  ∨ (x, n) = (3, 4) ∨ (x, n) = (-3, 4)
  ∨ (x, n) = (5, 5) ∨ (x, n) = (-5, 5)
  ∨ (x, n) = (11, 7) ∨ (x, n) = (-11, 7)
  ∨ (x, n) = (181, 15) ∨ (x, n) = (-181, 15)
\end{lstlisting}

\begin{figure}[!htbp]
\centering
\resizebox{\textwidth}{!}{%
\begin{tikzpicture}[
    >=Stealth,
    every node/.style={font=\small},
    result/.style={draw, rounded corners, fill=green!10,
                   text width=9em, align=center, inner sep=4pt},
    file/.style={draw, dashed, rounded corners, inner sep=10pt,
                 fill=blue!3},
    arr/.style={->, thick},
  ]
  \node[result] (rtheta) at (0,0)
    {\texttt{R},\,$\theta$,\,$\theta'$\\{\scriptsize $\theta^2{=}\theta{-}2$,\ $\theta\theta'{=}2$}};
  \node[result] (norm) at (0,2.1)
    {\texttt{norm\_eq}\\{\scriptsize $N(x{+}y\theta){=}x^2{+}xy{+}2y^2$}};

  \node[result] (ed) at (-3.4,4.4)
    {\texttt{instEuclidean\\Domain}\\{\scriptsize $R$ is Euclidean}};
  \node[result] (pid) at (-3.4,6.7)
    {\texttt{instPrincipal\\IdealRing}\\{\scriptsize $R$ is a PID}};
  \node[result] (ufd) at (-3.4,9.0)
    {\texttt{instUnique\\Factorization\\Monoid}\\{\scriptsize $R$ is a UFD}};

  \node[result] (units) at (3.4,4.4)
    {\texttt{units\_pm\_one}\\{\scriptsize $R^\times{=}\{\pm1\}$}};
  \node[result] (prime) at (3.4,9.0)
    {\texttt{theta\_prime}\\{\scriptsize $\theta,\theta'$ irreducible \& prime}};

  \draw[arr] (rtheta) -- (norm);
  \draw[arr] (norm) -- (ed);
  \draw[arr] (ed) -- (pid);
  \draw[arr] (pid) -- (ufd);
  \draw[arr] (norm) -- (units);
  \draw[arr] (norm) -- (prime);
  \draw[arr] (ufd) -- (prime);

  \node[result] (main) at (0,12.0)
    {\texttt{main\_m\_condition}\\{\scriptsize $-2\theta{+}1{=}\theta^m{-}\theta'^m$}};
  \draw[arr] (ufd) -- (main);
  \draw[arr] (units) -- (main);
  \draw[arr] (prime) -- (main);

  \node[result] (mod42) at (-3.4,14.6)
    {\texttt{odd\_case\_only\_\\three\_values\_mod\_42}\\{\scriptsize $m\equiv3,5,13\!\pmod{42}$}};
  \node[result] (unique) at (3.4,14.6)
    {\texttt{at\_most\_one\_\\m\_per\_class}\\{\scriptsize uniqueness per class}};
  \draw[arr] (main) -- (mod42);
  \draw[arr] (main) -- (unique);

  \node[result] (even) at (8.8,14.6)
    {\texttt{ramanujan\_nagell\_\\even\_pow\_factors}\\{\scriptsize even case: $n{=}4$}};

  \node[result, fill=green!25, text width=10em] (thm) at (0,17.4)
    {\textbf{\texttt{RamanujanNagell}}\\{\scriptsize the main theorem}};
  \draw[arr] (mod42) -- (thm);
  \draw[arr] (unique) -- (thm);
  \draw[arr] (even) -- (thm);

  \begin{scope}[on background layer]
    \node[file, fit=(rtheta)(norm)(ed)(pid)(ufd)(units)(prime),
          label={[font=\scriptsize\sffamily]left:Helpers.lean}] {};
    \node[file, fit=(main)(mod42)(unique)(even)(thm),
          label={[font=\scriptsize\sffamily]left:Basic.lean}] {};
  \end{scope}
\end{tikzpicture}%
}
\caption{Simplified dependency graph of the formalization. Boxes indicate
Lean declarations; dashed outlines group them by source file. Arrows point
from dependencies to dependents. The full interactive graph is available
in the \href{https://barindersbanwait.com/ramanujan_nagell/blueprint/dep_graph_document.html}{project blueprint}.}
\label{fig:dep-graph}
\end{figure}
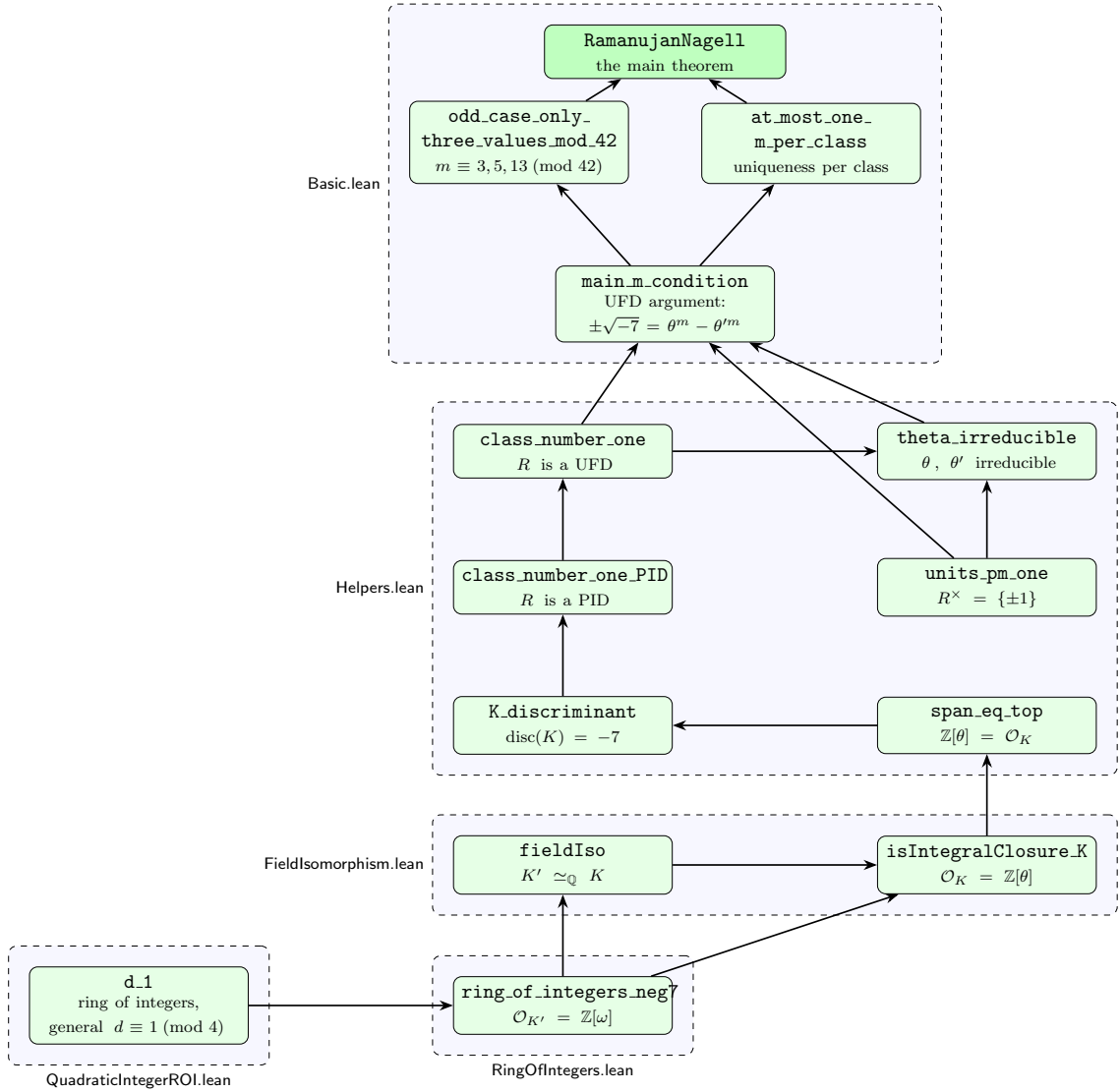

The formalization comprises two files, \texttt{Helpers.lean} and
\texttt{Basic.lean}, built directly on \Mathlib{}.
\Cref{fig:dep-graph} shows the dependency graph of the principal
results, simplified from the full blueprint graph (which contains
approximately 70~nodes).
An interactive version is available in the project's online blueprint.\footnote{\url{https://barindersbanwait.com/ramanujan_nagell/blueprint/dep_graph_document.html}}

The key results along this chain, roughly in dependency order, are:
\begin{enumerate}[label=(\arabic*)]
  \item \textbf{The ring $R$ and its arithmetic}
    (\texttt{Helpers.lean}).
    We work in $R = \Z[\tfrac{1+\sqrt{-7}}{2}]$, realized as the
    \Mathlib{} type \lean{QuadraticAlgebra}\,$\Z$\,$(-2)$\,$1$, with
    $\theta = \langle 0,1\rangle$ and $\theta' = \langle 1,-1\rangle$.
    The identities $\theta^2 = \theta - 2$, $\theta\theta' = 2$,
    $\theta + \theta' = 1$ hold by \texttt{rfl}, and the norm is
    $N(x+y\theta) = x^2 + xy + 2y^2$.

  \item \textbf{$R$ is a Euclidean domain, hence a UFD}
    (\texttt{Helpers.lean}).
    With $N$ as Euclidean function, $R$ is a Euclidean domain
    (\lean{instEuclideanDomain}), hence a principal ideal domain
    (\lean{instPrincipalIdealRing}) and a unique factorization domain
    (\lean{instUniqueFactorizationMonoid}).

  \item \textbf{Units and primes}
    (\texttt{Helpers.lean}).
    The only units of $R$ are $\pm 1$ (\lean{units_pm_one}), by
    completing the square in $4N = (2x+y)^2 + 7y^2$.
    Since $N(\theta) = N(\theta') = 2$ is prime, $\theta$ and $\theta'$
    are irreducible (\lean{theta_irreducible}), hence prime
    (\lean{theta_prime}), and non-associated.

  \item \textbf{The factorization and sign argument}
    (\texttt{Basic.lean}).
    Factoring $(x^2+7)/4 = \theta^m\theta'^m$ in the UFD~$R$ and pinning
    down the sign yields the key identity
    $-2\theta+1 = \theta^m - \theta'^m$ (\lean{main_m_condition}).

  \item \textbf{The mod-$42$ reduction and uniqueness argument}
    (\texttt{Basic.lean}).
    The binomial expansion gives $-2^{m-1} \equiv m \pmod{7}$, reducing
    the exponent to three residue classes modulo~$42$
    (\lean{odd_case_only_three_values_mod_42}); a $7$-adic valuation
    argument then shows each class contributes at most one solution
    (\lean{at_most_one_m_per_class}). The even case is elementary.
\end{enumerate}
These results fall into two categories.
Items (1)--(3) constitute the \emph{algebraic infrastructure}
(\texttt{Helpers.lean}): the properties of~$R$ that a textbook takes for
granted or dispatches in a line, the key results being
\lean{instEuclideanDomain} and \lean{units_pm_one}.
Items (4) and~(5) form the \emph{proof proper} (\texttt{Basic.lean}):
the unique factorization argument, the mod-$42$ reduction, and the
uniqueness per residue class, the main results being
\lean{main_m_condition}, \lean{odd_case_only_three_values_mod_42},
and \lean{at_most_one_m_per_class} (the last of which relies on a
non-trivial $7$-adic valuation argument).
We highlight results from both categories in \Cref{sec:statements}
and discuss the most challenging aspects in \Cref{sec:challenges}.

\subsection{Module structure and statistics}\label{sec:modules}

\Cref{tab:modules} summarizes the contents of \texttt{Helpers.lean}
and \texttt{Basic.lean}.

\begin{table}[ht]
\centering
\caption{Per-module breakdown of the formalization. LOC = lines of code; Decls = number of Lean declarations (definitions, lemmas, theorems, and instances).}
\label{tab:modules}
\begin{tabular}{@{}l p{0.50\textwidth} r r@{}}
\toprule
File & Purpose & LOC & Decls \\
\midrule
\texttt{Helpers.lean}
  & Ring $R$ and norm; Euclidean domain $\Rightarrow$ PID $\Rightarrow$ UFD; units $\pm1$; irreducibility and primality of $\theta,\theta'$ & 506 & 42 \\[3pt]
\texttt{Basic.lean}
  & Factorization and sign argument, mod-$42$ reduction, uniqueness per class, assembly of the main theorem & 1\,305 & 41 \\
\midrule
\textbf{Total} & & \textbf{1\,811} & \textbf{83} \\
\bottomrule
\end{tabular}
\end{table}

\texttt{Helpers.lean} sets up the ring $R = \Z[\tfrac{1+\sqrt{-7}}{2}]$ as
the \Mathlib{} quadratic algebra \lean{QuadraticAlgebra}\,$\Z$\,$(-2)$\,$1$
and develops the algebraic properties the main argument needs: the norm
form, the Euclidean-domain structure (with the resulting PID and UFD
instances), the determination of the unit group $R^\times = \{\pm 1\}$, and
the irreducibility and primality of $\theta$ and~$\theta'$.

\texttt{Basic.lean} contains the proof itself: the unique factorization
argument (\lean{main_m_condition}), the binomial-expansion and
mod-$42$ reduction (\lean{odd_case_only_three_values_mod_42}),
the uniqueness per residue class (\lean{at_most_one_m_per_class}),
and the final assembly into the main theorem.

The project builds against \Mathlib{} commit \texttt{fa08db7d}
on Lean~\texttt{v4.30.0}.
The formalization contains no \texttt{sorry} placeholders; after
\texttt{lake build}, running
\begin{quote}
\begin{lstlisting}[basicstyle=\ttfamily\small]
echo 'import RamanujanNagell.Basic
#print axioms RamanujanNagell' | lake env lean --stdin
\end{lstlisting}
\end{quote}
confirms that the main theorem depends only on the three standard
Lean axioms (\texttt{propext}, \texttt{Classical.choice},
\texttt{Quot.sound}).

\subsection{Key formal statements}\label{sec:statements}

We present the formal statements of several key results, organized by
the two categories identified above.

\subsubsection*{Algebraic infrastructure}

\Mathlib{}'s \lean{QuadraticAlgebra}\,\texttt{A\,d\,e} is the free
\texttt{A}-module on $1$ and a generator~$\omega$ with
$\omega^2 = d + e\omega$.  Taking $\texttt{A} = \Z$, $d = -2$, $e = 1$
gives exactly $R = \Z[\tfrac{1+\sqrt{-7}}{2}]$, with
$\theta = \langle 0,1\rangle$ and $\theta' = \langle 1,-1\rangle$ as plain
elements --- no integral-closure detour and no subtype wrapping.  The basic
identities are then definitional:

\begin{lstlisting}
lemma theta_sq : θ ^ 2 = θ - 2          -- by rfl
lemma theta_mul_theta' : θ * θ' = 2     -- by rfl
\end{lstlisting}

Unique factorization is obtained explicitly.  With the norm
$N(x+y\theta) = x^2 + xy + 2y^2$ as Euclidean function, $R$ is a Euclidean
domain, from which \Mathlib{} derives the PID and UFD instances.

\begin{lstlisting}
noncomputable instance instEuclideanDomain : EuclideanDomain R
instance instPrincipalIdealRing : IsPrincipalIdealRing R
instance instUniqueFactorizationMonoid : UniqueFactorizationMonoid R
\end{lstlisting}

The division algorithm rounds $a\bar{b}/N(b)$ to a nearby element of $R$;
the bound ensuring the remainder has strictly smaller norm is a short
computation (\Cref{sec:euclidean}).

That $R^\times = \{\pm 1\}$ follows by completing the square: a unit has
norm~$1$, and $4 = (2x+y)^2 + 7y^2$ forces $y = 0$ and $x = \pm 1$.

\begin{lstlisting}
lemma units_pm_one (u : Rˣ) : u = 1 ∨ u = -1
\end{lstlisting}

Finally $N(\theta) = N(\theta') = 2$ is prime in $\Z$, so $\theta$ and
$\theta'$ are irreducible, hence prime in the UFD~$R$, and non-associated
--- completing the algebraic setup for the main argument.

\begin{lstlisting}
lemma theta_irreducible : Irreducible θ
lemma theta_prime : Prime θ
lemma theta_theta'_not_associated : ¬ Associated θ θ'
\end{lstlisting}

\subsubsection*{The proof proper}

The core step factors $(x^2+7)/4 = 2^m$ in the UFD~$R$
and uses unique factorization plus the coprimality of $\theta$ and
$\theta'$ to conclude the key identity $-2\theta+1 = \theta^m - \theta'^m$.

\begin{lstlisting}
theorem main_m_condition :
  ∀ x : ℤ, ∀ m : ℕ, Odd m → m ≥ 3 →
    (x ^ 2 + 7) / 4 = 2 ^ m →
    (-2*θ + 1 = θ^m - θ'^m)
\end{lstlisting}

A binomial expansion and reduction modulo~$7$, combined with Fermat's
little theorem ($2^6 \equiv 1 \pmod{7}$), restricts the exponent to
three residue classes.

\begin{lstlisting}
theorem odd_case_only_three_values_mod_42 :
  ∀ x : ℤ, ∀ n : ℕ, Odd n → n ≥ 5 →
    x ^ 2 + 7 = 2 ^ n →
    (n-2) % 42 = 3 ∨ (n-2) % 42 = 5 ∨ (n-2) % 42 = 13
\end{lstlisting}

Finally, the uniqueness argument shows that each residue class
contributes at most one solution (the formal counterpart of
\Cref{sec:uniqueness-math}).  The proof proceeds by contradiction:
if $m_1 \equiv m_2 \pmod{42}$ both satisfy the theta equation, then
setting $d = m_2 - m_1$ yields a $7$-adic valuation contradiction via
the identity $P \cdot B_d = 1 - 7A'_d - 2^d$, where $P$ is
the trace $\theta^{m_1} + \theta'^{m_1}$ (not divisible by~$7$)
and $v_7(B_d) = v_7(d)$ by a binomial coefficient analysis.

\begin{lstlisting}
lemma at_most_one_m_per_class (m₁ m₂ : ℕ)
    (h₁_odd : Odd m₁) (h₂_odd : Odd m₂)
    (h₁_ge : m₁ ≥ 3) (h₂_ge : m₂ ≥ 3)
    (h_cong : m₁ % 42 = m₂ % 42)
    (h₁_theta : -2 * θ + 1 = θ ^ m₁ - θ' ^ m₁)
    (h₂_theta : -2 * θ + 1 = θ ^ m₂ - θ' ^ m₂) :
    m₁ = m₂
\end{lstlisting}

\subsection{An earlier approach via the ring of integers}\label{sec:earlier}

The proof formalized here is a simplification of an earlier
version of the project, preserved on the \texttt{legacy} branch of the
repository,\footnote{Reachable from the project page,
\url{https://barindersbanwait.com/ramanujan_nagell/}; the snapshot is
tagged \texttt{legacy-OK-formalization}.} which proved the same theorem
by a different route. That version worked in the number field
$K = \Q(\sqrt{-7})$, presented as
\lean{QuadraticAlgebra}\,$\Q$\,$(-2)$\,$1$; it
identified the ring of integers $\cO_K$ using the integral-closure
criterion for quadratic fields formalized by Brasca and
Monticone~\cite{brasca-monticone-quadratic}, proved $\cO_K$ a principal
ideal domain from the Minkowski bound on its discriminant (via
\Mathlib{}'s \lean{isPrincipalIdealRing_of_abs_discr_lt}), and determined
the unit group through Dirichlet's unit theorem. Reconciling the two
natural presentations of $\Q(\sqrt{-7})$---the one generated by
$\sqrt{-7}$ used by Brasca and Monticone, and the one generated by
$\tfrac{1+\sqrt{-7}}{2}$ that gives the ring of integers---further
required an explicit $\Q$-algebra isomorphism transporting the
integral-closure property between them.

That development ran to roughly $3{,}570$ lines\footnote{although no serious attempt was made to shorten this}, and the comparison with
the present proof is instructive. On paper, citing that $\Q(\sqrt{-7})$
has class number one is the shorter route, and proving the norm-Euclidean property by hand looks like extra work. In formalization the balance inverts:
\Mathlib{} yields class number one only through the Minkowski bound,
which presupposes the entire number-field apparatus---the ring of
integers, its discriminant, and, for the units, Dirichlet's unit
theorem---together with the isomorphism transport between the two
presentations of the field. The norm-Euclidean route needs none of this:
it is a self-contained division algorithm of about $170$ lines
(\Cref{sec:euclidean}), and the units come from an elementary
completing-the-square argument (\Cref{sec:units}). The present proof is
consequently about half the size, with the heaviest abstract
infrastructure replaced by elementary computation.

We nonetheless retain the earlier formalization, because its
infrastructure---the identification of the ring of integers, the
discriminant, and the principal-ideal-domain criterion---applies to
every imaginary quadratic field, and so is the natural foundation for
generalizing the method to other equations of Ramanujan--Nagell type,
where the norm-Euclidean shortcut is unavailable.

\section{Challenges in formalization}\label{sec:challenges}

We now discuss the aspects of the formalization that were unexpectedly
difficult or that diverge most sharply from the textbook presentation.
The largest single difficulty of an earlier version of this
project---transporting the ring of integers between two presentations
of $\Q(\sqrt{-7})$---was \emph{removed} by working directly in
$R = \texttt{QuadraticAlgebra}\;\Z\,(-2)\,1$, where $\theta$, $\theta'$,
and their defining relations are available with no integral-closure
detour. What remains is dominated instead by a single task that a
textbook discharges in one word: making $R$ a Euclidean domain.

\subsection{Constructing the Euclidean domain}\label{sec:euclidean}

\Cref{sec:ufd} establishes that $R$ is norm-Euclidean by a rounding
argument in $\Q(\sqrt{-7})$. Formalizing that argument is the largest
single component of the algebraic infrastructure (\texttt{Helpers.lean})---and
the part forward-referenced from \Cref{sec:statements}---since \Mathlib{} neither
knows the fact nor provides any criterion for the norm-Euclidean property: we must
exhibit the division algorithm explicitly, producing for $a, b \in R$
with $b \neq 0$ an element $q \in R$ with $N(a - bq) < N(b)$.

Here the formalization deviates from the textbook argument in one
respect. The exact quotient $a/b$ lives in $\Q(\sqrt{-7})$, not in $R$;
rather than construct elements of the ambient field in \Lean{}, we clear
the denominator first. Multiplying by the conjugate $\bar b$ and writing
$s = a\bar b$, multiplicativity of the norm gives
\[
  N(b)\,N(a - bq) \;=\; N\!\bigl(s - N(b)\,q\bigr)
\]
(\lean{N_mul_norm_rem_eq}), which reduces the problem to rounding the
\emph{element} $s$ to the nearest multiple of the rational integer
$N(b)$. The quotient is then built by rounding the $\sqrt{-7}$-coordinate
of $s / N(b)$ first, and re-rounding the rational coordinate after
shifting by half the first residual:

\begin{lstlisting}
noncomputable def quot (a b : R) : R :=
  let N := QuadraticAlgebra.norm b
  if N = 0 then 0 else
    let s := a * star b              -- s = a * conj(b)
    let n := round ((s.im : ℚ) / N)  -- round the imaginary part first
    let m := round ((2*s.re + s.im - N*n) / (2*N))
    ⟨m, n⟩
\end{lstlisting}

Writing $u$ and $v$ for the two integer residuals of $s$ after this
rounding, the heart of the matter (\lean{sixteen_norm_rem_le},
about 65~lines) is the algebraic identity
\[
  16\,(u^2 + uv + 2v^2) \;=\; 4\,(2u+v)^2 + 7\,(2v)^2
\]
together with the rounding bounds $(2u+v)^2 \le N(b)^2$ and
$(2v)^2 \le N(b)^2$, each coming from
$|\!\operatorname{round}(x) - x| \le \tfrac12$. Since
$N(b)\,N(a-bq) = u^2 + uv + 2v^2$, these combine to
\[
  16\,N(b)\,N(a-bq) \le 11\,N(b)^2,
  \qquad\text{hence}\qquad
  N(a-bq) \le \tfrac{11}{16}\,N(b) < N(b).
\]
The constant $\tfrac{11}{16}$ is precisely the one from the rounding
argument of \Cref{sec:ufd}, here recovered entirely within~$\Z$. The
Euclidean function itself is $|N|$, valued in $\N$ (\lean{normMeasure}).
With \lean{quot}, the remainder \lean{rem}, and this bound, the three
instances follow, \Mathlib{} supplying the last two from the first:

\begin{lstlisting}
noncomputable instance instEuclideanDomain : EuclideanDomain R
instance instPrincipalIdealRing : IsPrincipalIdealRing R
instance instUniqueFactorizationMonoid : UniqueFactorizationMonoid R
\end{lstlisting}

The whole development runs to roughly 170~lines---about a third of
\texttt{Helpers.lean}---and is the formal cost of the textbook's single
word ``norm-Euclidean''. The instance is necessarily
\texttt{noncomputable}, since \lean{quot} rounds in~$\Q$; this costs
nothing here, as the structure serves only to produce the PID and UFD
\emph{instances}, never to compute a quotient.

\subsection{Determining the unit group}\label{sec:units}

The determination of $R^\times = \{\pm 1\}$ is the sharpest
illustration of how the choice of setting reshapes a formalization. In the earlier $\cO_K$-based version of this project it was the single most
labour-intensive result---some 236~lines---because the natural argument
runs through Dirichlet's unit theorem (every unit of an imaginary
quadratic field is a root of unity), a totient bound
$\varphi(n) \le [K:\Q] = 2$ restricting the order to
$n \in \{1,2,3,4,6\}$, and then the elimination of the orders $3$, $4$,
and~$6$ via a Cayley--Hamilton degree bound and a case analysis on
minimal polynomials. Dirichlet's theorem is itself a central result of
an undergraduate course in algebraic number theory, and although it is
available in \Mathlib{}, the totient bound and the root-of-unity
elimination are not, and reconstructing them dominated the file.

Working in $R$ with the explicit norm form dispenses with all of this. A
unit $u$ has $N(u) \in \Z^\times = \{\pm 1\}$, and since $N$ is
nonnegative $N(u) = 1$; completing the square in
$4N(u) = (2x + y)^2 + 7y^2 = 4$ (where $u = x + y\theta$) forces $y = 0$,
then $x^2 = 1$, so $u = \pm 1$. The whole proof is under 30~lines, the
only non-routine step being a single \texttt{nlinarith} call---given the
hints $\texttt{sq\_nonneg}\;y$ and $\texttt{sq\_nonneg}\,(2x+y)$---to
conclude $y = 0$ from the bound.

\begin{lstlisting}
lemma units_pm_one (u : Rˣ) : u = 1 ∨ u = -1
\end{lstlisting}

\noindent What had been the most demanding result of the formalization
became one of the easiest.

\subsection{Arithmetic in $R$}\label{sec:arithmetic}

Working in $R$ rather than~$\Z$ is mostly invisible, but it occasionally
intrudes. On the positive side, the defining relations are now
definitional: $\theta^2 = \theta - 2$, $\theta\theta' = 2$, and
$\theta + \theta' = 1$ all hold by \texttt{rfl}, where the earlier
$\cO_K$ presentation needed a multi-step coercion through~$K$ for each.

The residual friction is that the integer tactics \texttt{ring},
\texttt{linarith}, and \texttt{omega} operate on~$\Z$, not on~$R$, so an
integer (in)equality concealed inside an $R$-equation must be extracted
by hand. The standard move is to apply \lean{congrArg} with
\lean{QuadraticAlgebra.re} or \lean{QuadraticAlgebra.im} (or
\lean{QuadraticAlgebra.ext_iff} to split both coordinates at once),
reducing an $R$-statement to its two integer components, which
\texttt{simp} and \texttt{omega} then finish. For instance
\lean{factor_not_unit_left} (25~lines)---the formal content of the
textbook's ``the imaginary part is nonzero''---uses \lean{units_pm_one}
to reduce a putative unit factor to $\pm 1$, then compares imaginary
parts via \lean{congrArg QuadraticAlgebra.im} to reach a contradiction,
since $2\theta - 1$ has nonzero imaginary part while $\pm 2^m$ does not.

One instance-level wrinkle is worth recording. The layered structure of
\lean{QuadraticAlgebra} (a quadratic algebra over~$\Z$, itself a
commutative ring, carrying a norm and a star operation) gives rise to
instance-resolution chains deep enough that the project raises
\texttt{maxSynthPendingDepth} to~$3$ in its \texttt{lakefile} so that
typeclass inference will pursue them. This is the only typeclass-inference
adjustment the formalization requires---in marked contrast to the
$\cO_K$ development, where reconciling competing
\texttt{Algebra}\,$\Q$\,$K$ instances demanded explicit transports
throughout.

\subsection{The binomial expansion and mod-42 reduction}\label{sec:mod42}

The reduction from an arbitrary odd exponent~$m$ to three residue
classes modulo~$42$ falls into two steps, both now carried out directly
in~$R$. Introduce $\delta = 2\theta - 1$, which satisfies
$\delta^2 = -7$ (a three-line \texttt{calc} from $\theta^2 = \theta - 2$);
since the key identity gives $\theta^m - \theta'^m = -\delta$,
multiplying by $2^m$ and using $\theta' = 1 - \theta$ yields
$(2\theta)^m - (2(1-\theta))^m = -2^m\delta$.

The first step (\lean{expand_by_binomial}, about 140~lines) applies
the binomial theorem to $(2\theta)^m$ and $(2(1-\theta))^m$ via
\lean{add_pow}, subtracts, and separates the result into even- and
odd-indexed terms. The even terms cancel (using \lean{Even.neg_pow}) and
the odd terms double (\lean{Odd.neg_pow}); in \Lean{} the separation
uses \lean{Finset.sum_filter_add_sum_filter_not} to partition the index
set, then a \lean{Finset.sum_bij} to reindex the surviving odd indices
$k = 2j+1$ to consecutive indices~$j$. The \lean{sum_bij} call alone
discharges four obligations (membership, injectivity, surjectivity, and
the per-index value identity, the last using $\delta^2 = -7$). Reducing
the resulting integer identity modulo~$7$ gives
$-2^{m-1} \equiv m \pmod{7}$.

The second step (\lean{odd_case_only_three_values_mod_42}, about
50~lines) combines this congruence with Fermat's little theorem in the
form $2^6 \equiv 1 \pmod{7}$---established by a short induction on the
exponent---and the constraint that $m$ is odd. The proof splits on
$m \bmod 6 \in \{1, 3, 5\}$, and in each case \texttt{omega} performs the
CRT-style combination of the mod-$6$ and mod-$7$ data into the mod-$42$
conclusion $m \equiv 3, 5, 13 \pmod{42}$.

\subsection{The uniqueness argument and $7$-adic valuations}\label{sec:uniqueness}

The longest single proof in the formalization is
\lean{at_most_one_m_per_class} (about 225~lines), which shows that each
residue class modulo~$42$ contributes at most one solution, formalizing
the argument of \Cref{sec:uniqueness-math}. The proof proceeds by
contradiction: assuming $m_1 < m_2$ with $m_1 \equiv m_2 \pmod{42}$, one
sets $d = m_2 - m_1$ and $l = v_7(d) \ge 1$, and derives a $7$-adic
valuation contradiction.

Its core is the algebraic identity
\[
  P \cdot B_d = 1 - 7 A'_d - 2^d,
\]
where $P = \theta^{m_1} + \theta'^{m_1}$ is the trace (an integer
coprime to~$7$), $B_d = \sum_j \binom{d}{2j+1}(-7)^j$ is the odd-indexed
binomial sum, with $v_7(B_d) = l$, and $A'_d$ is the corresponding
even-indexed sum. This identity is built in~$R$ by the method of \Cref{sec:mod42}:
introduce $\delta = 2\theta - 1$ (so $\delta^2 = -7$), expand
$(2\theta)^d$ and $(2(1-\theta))^d$ with \lean{add_pow}, and separate the
odd-indexed terms from the even-indexed ones, reindexing the odd half by
\lean{sum_bij} and the even half by \lean{Finset.sum_image}. The
resulting equation in~$R$ is then transferred to~$\Z$: since the cast
$\Z \to R$ is injective, it suffices to compare \lean{re} components,
which \lean{congrArg} extracts. Constructing and lifting this identity
accounts for roughly 145 of the proof's lines.

The contradiction then follows from three valuation inputs:
\lean{lemma_A_binomial_valuation} gives $v_7(B_d) = l$;
\lean{even_binomial_valuation} gives $7^l \mid A'_d$; and the
lifting-the-exponent lemma---which evaluates $v_7(a^k - b^k)$ as
$v_7(a - b) + v_7(k)$ (\Cref{sec:uniqueness-math}), and is available in
\Mathlib{} as \lean{Int.emultiplicity_pow_sub_pow}---applied to
$8^{d/3} - 1$ since $3 \mid d$, gives $7^{l+1} \mid 2^d - 1$. Substituted
into the identity,
these show $7^{l+1} \mid P \cdot B_d$; as $\gcd(P, 7) = 1$, this forces
$7^{l+1} \mid B_d$, contradicting $v_7(B_d) = l$.

The trace coprimality $7 \nmid P$ is itself a small lemma
(\lean{trace_seq_not_dvd_seven}): the integers
$P_k = \theta^k + \theta'^k$ satisfy the recurrence
$P_{k+2} = P_{k+1} - 2P_k$ with $P_0 = 2$, $P_1 = 1$, and modulo~$7$ this
sequence is periodic of period~$3$ with values $\{2, 1, 4\}$, none of
which is~$0$.

\section{Lessons for formalizing algebraic number theory}\label{sec:lessons}

\subsection{What \Mathlib{} provided}\label{sec:mathlib-provided}

Several components of the proof relied heavily on infrastructure
already present in \Mathlib{}:
\begin{itemize}
  \item \emph{$p$-adic valuations.}  The \lean{padicValNat} and
    \lean{multiplicity} APIs, together with the Lifting-the-Exponent
    lemma (\lean{Int.emultiplicity_pow_sub_pow}), were essential
    for the uniqueness argument in \cref{sec:uniqueness}.
  \item \emph{The factorization hierarchy.}  Once we supplied the
    Euclidean-domain structure on~$R$ (\cref{sec:euclidean}),
    \Mathlib{} furnished the rest of the chain for free: the
    implications \lean{EuclideanDomain} $\Rightarrow$
    \lean{IsPrincipalIdealRing} $\Rightarrow$
    \lean{UniqueFactorizationMonoid}, and the prime theory of a UFD that
    promotes $N(\theta) = N(\theta') = 2$ being prime in~$\Z$ to
    primality of $\theta$ and~$\theta'$ in~$R$.
  \item \emph{Quadratic extensions.}  \Mathlib{}'s \lean{QuadraticAlgebra}
    provided a ready-made type for $R = \Z[\tfrac{1+\sqrt{-7}}{2}]$,
    including the norm form (via its \lean{NormDeterminant} API), the
    star/conjugation operation, and the coordinate projections on which
    most of the elementary arithmetic rests.
  \item \emph{Finset combinatorics.}  The binomial theorem
    (\lean{add_pow}), \lean{Finset.sum_filter_add_sum_filter_not}, and
    \lean{Finset.sum_bij} were used extensively in the binomial
    expansion and reindexing arguments of \cref{sec:mod42}.
\end{itemize}

\subsection{What was missing or painful}\label{sec:missing}

Despite the breadth of \Mathlib{}, two gaps required significant
manual effort:
\begin{itemize}
  \item \emph{The Norm-Euclidean property.}  \Mathlib{} has no criterion or
    decision procedure for whether a quadratic ring is norm-Euclidean,
    so the division algorithm and the remainder bound underlying the
    Euclidean-domain instance had to be built by hand
    (\cref{sec:euclidean}). This is the largest body of work in the
    project with no counterpart in a textbook proof, which simply
    asserts that $\Q(\sqrt{-7})$ is one of the norm-Euclidean imaginary
    quadratic fields.
  \item \emph{Arithmetic in $R$.}  The absence of decision procedures
    for quadratic rings means that identities trivial over~$\Z$ (where
    \texttt{ring} and \texttt{omega} suffice) require extracting integer
    components by hand when carried out in~$R$ (see
    \cref{sec:arithmetic}).
\end{itemize}

\subsection{AI-assisted formalization}\label{sec:ai}

A distinctive feature of this project is the significant role played
by generative AI tools throughout the formalization process.  We made
extensive use of two systems: \emph{Claude Code}~\cite{ClaudeCode},
Anthropic's command-line coding agent equipped with Lean~4-specific
skills developed by Cameron Freer~\cite{FreerLean4Skills}, and
\emph{Aristotle}~\cite{Aristotle}, an AI system specialized for
Lean~4 proof synthesis.

These tools assisted at multiple levels of the formalization:
\begin{itemize}
  \item \emph{Tactic suggestion.}  Given a proof state, both tools
    could frequently suggest the correct next tactic or a short tactic
    sequence, particularly for routine goals involving \texttt{ring},
    \texttt{omega}, \texttt{simp}, and \texttt{norm\_num}.  This
    eliminated much of the trial-and-error overhead that typically
    dominates interactive theorem proving.
  \item \emph{Proof repair.}  When \Mathlib{} updates introduced
    breaking changes to API names or signatures, Claude Code could
    read compiler errors and automatically update proofs to use the
    new API, a task that would otherwise require manually searching
    \Mathlib{} documentation.
  \item \emph{Lemma discovery.}  Searching \Mathlib{}'s large API
    surface for the right lemma is a well-known bottleneck in Lean
    formalization.  The AI tools could often identify relevant lemmas
    from natural-language descriptions of the desired result,
    significantly reducing search time.
  \item \emph{Boilerplate generation.}  Structurally repetitive proof
    steps --- such as the per-index obligations of the \lean{sum_bij}
    reindexings (\cref{sec:mod42}) or the parallel case analyses in the
    mod-$42$ reduction and uniqueness arguments --- were generated much
    faster with AI assistance, even when each case required minor manual
    corrections.
\end{itemize}

We emphasize that every proof in the formalization is ultimately
verified by Lean's kernel; the AI tools accelerate the \emph{writing}
of proofs but do not affect their \emph{correctness}.  In our
experience, these tools massively lowered the barrier to entry for
Lean formalization, and we expect
AI-assisted formalization to become standard practice as these tools mature.

\subsection*{Acknowledgments}

We are especially grateful to Michael Stoll for helpful comments on both the
exposition and the code, and for insisting that things be done with $R$ only, which
massively simplified and shortened the formalization. We are very grateful to Riccardo
Brasca for code and guidance supplied on the Lean Zulip. We thank Andrew O'Desky,
Ashvni Narayanan, and Mar\'ia In\'es de Frutos Fern\'andez for helpful discussions. We
are also grateful to Xinze Bryan Li for some helpful contributions to the formalization.

We acknowledge the use of Claude
Code~\cite{ClaudeCode} with Lean~4 skills by Cameron Freer~\cite{FreerLean4Skills}, and Aristotle~\cite{Aristotle}, for AI-assisted proof development.  This paper was written while the author was at the Banff International Research Station (BIRS) for the workshop \emph{DANGER: Data, Numbers, and Geometry}. We thank the workshop organizers, Edward Hirst and Elli Heyes, for the invitation, and BIRS for the excellent working conditions.


\renewcommand{\refname}{
\begin{center}
\normalsize \bf References
\end{center}
}

\bibliographystyle{plainnat}
\bibliography{references}

\end{document}